\documentclass[12pt]{amsart}
\usepackage{amsmath,amssymb,amsfonts,amsthm,amscd,indentfirst}
\usepackage{hyperref}\usepackage{xcolor}

\numberwithin{equation}{section}

\newtheorem{prop}{Proposition}[section]
\newtheorem{theo}[prop]{Theorem}
\newtheorem{lemm}[prop]{Lemma}

\def\and{\quad{\rm and}\quad}

\def\<{\langle}
\def\>{\rangle}

\usepackage{graphicx}

\begin{document}
\title{On the asymptotic Plateau problem in hyperbolic space}

\author[Siyuan Lu]{Siyuan Lu}
\address{Department of Mathematics and Statistics, McMaster University, 
1280 Main Street West, Hamilton, ON, L8S 4K1, Canada.}
\email{siyuan.lu@mcmaster.ca}
\thanks{The author was supported in part by NSERC Discovery Grant.}

\begin{abstract}
In this paper, we solve the asymptotic Plateau problem in hyperbolic space for constant $\sigma_{n-1}$ curvature, i.e. the existence of a complete hypersurface in $\mathbb{H}^{n+1}$ satisfying $\sigma_{n-1}(\kappa)=\sigma\in (0,n)$ with a prescribed asymptotic boundary $\Gamma$. The key ingredient is the curvature estimates. Previously, this is only known for $\sigma_0<\sigma<n$, where $\sigma_0$ is a positive constant.
\end{abstract}

\maketitle

\section{Introduction}

Let $\mathbb{H}^{n+1}$ be the hyperbolic space and let $\partial_\infty\mathbb{H}^{n+1}$ be the ideal boundary of $\mathbb{H}^{n+1}$ at infinity. The asymptotic Plateau problem in hyperbolic space asks to find a complete hypersurface of constant curvature in $\mathbb{H}^{n+1}$ with prescribed asymptotic boundary at infinity. More precisely, given a closed embedded smooth $(n-1)$-dimensional submanifold $\Gamma\subset \partial_\infty\mathbb{H}^{n+1}$, one seeks a complete hypersurface $\Sigma$ in $\mathbb{H}^{n+1}$ satisfying
\begin{align*}
f(\kappa)=\sigma,\quad \partial\Sigma=\Gamma,
\end{align*}
where $f$ is a smooth symmetric function of $n$ variables, $\kappa=(\kappa_1,\cdots,\kappa_n)$ are the principal curvatures of $\Sigma$ and $\sigma$ is a constant.

The asymptotic Plateau problem was first studied by Anderson \cite{A82,A83} and Hardt and Lin \cite{HL} for area-minimizing varieties using geometric measure theory. Their results were extended by Tonegawa \cite{T} to hypersurfaces of constant mean curvature. The asymptotic Plateau problem for constant mean curvature was studied by Lin \cite{Lin}, Nelli and Spruck \cite{NS} and Guan and Spruck \cite{GS00} using PDE methods. The asymptotic Plateau problem for constant Gauss curvature was studied by Labourie \cite{L} in $\mathbb{H}^3$ and by Rosenberg and Spruck \cite{RS} in $\mathbb{H}^{n+1}$. For a broad class of $f$ defined in the positive cone $\Gamma_n$, the asymptotic Plateau problem was completely solved via works of Guan, Spruck, Szapiel and Xiao \cite{GSS, GS11, GSX}. For general $f$ satisfying natural structure conditions, the asymptotic Plateau problem was solved by Guan and Spruck \cite{GS10} if $\sigma>\sigma_0$, where $\sigma_0$ is a positive constant.

After the work of Guan and Spruck \cite{GS10}, it is of great interest to study the remaining case, i.e. $0<\sigma\leq \sigma_0$. As pointed out in \cite{GS10}, the only missing piece is the curvature estimates. In a recent work by Wang \cite{Wang}, he was able to obtain the curvature estimates for $f=\frac{\sigma_k}{\sigma_{k-1}}$ following the method in \cite{GS10}, where $\sigma_k$ denotes the $k$-th elementary symmetric function.  In contrast, the case for $f=\sigma_k$ remains open.

\medskip

In this paper, we study the asymptotic Plateau problem for an important case $f=\sigma_{n-1}$. For $n=2$, it reduces to the mean curvature case, which was completely solved by previous works. Therefore we will restrict ourselves to the case $n\geq 3$.  

Before we state our main theorems, let us introduce some notations. Let 
\begin{align*}
\mathbb{H}^{n+1}=\{(x,x_{n+1})\in \mathbb{R}^{n+1}:x_{n+1}>0\},\quad ds^2=\frac{1}{x_{n+1}^2}\sum_{i=1}^{n+1}dx_i^2,
\end{align*}
denote the upper half-space model of $\mathbb{H}^{n+1}$. Then $\partial_\infty\mathbb{H}^{n+1}$ is naturally identified with $\mathbb{R}^n=\mathbb{R}^n\times \{0\}\subset \mathbb{R}^{n+1}$. 

Denote Garding's $\Gamma_k$ cone by
\begin{align*}
\Gamma_k=\{\lambda\in \mathbb{R}^n:\sigma_j(\lambda)>0,1\leq j\leq k\}.
\end{align*}

We now state our main theorems.
\begin{theo}\label{theo-estimate}
For $n\geq 3$, let $\Gamma=\partial\Omega$, where $\Omega$ is a bounded smooth domain in $\mathbb{R}^n$ with nonnegative mean curvature and $\sigma\in (0,n)$. Suppose $\Sigma$ is a $C^4$ vertical graph over $\Omega$ in $\mathbb{H}^{n+1}$ satisfying
\begin{align}\label{sigma_n-1}
\sigma_{n-1}(\kappa)=\sigma,\quad \partial\Sigma=\Gamma,
\end{align}
where $\kappa\in \Gamma_{n-1}$. Then we have
\begin{align*}
\max_{x\in \Sigma, 1\leq i\leq n}|\kappa_i(x)|\leq C+C\max_{x\in \partial\Sigma, 1\leq i\leq n}|\kappa_i(x)|,
\end{align*}
where $C$ is a constant depending only on $n, \Omega$ and $\sigma$.
\end{theo}

As an application, we solve the asymptotic Plateau problem for constant $\sigma_{n-1}$ curvature for all $\sigma$.
\begin{theo}\label{theo-exist}
For $n\geq 3$, let $\Gamma=\partial\Omega$, where $\Omega$ is a bounded smooth domain in $\mathbb{R}^n$ with nonnegative mean curvature and $\sigma\in (0,n)$. Then there exists a complete hypersurface $\Sigma$ in $\mathbb{H}^{n+1}$ satisfying
\begin{align*}
\sigma_{n-1}(\kappa)=\sigma,\quad \partial\Sigma=\Gamma.
\end{align*}
\end{theo}

To obtain the curvature estimates, we adopt a new test function and make full use of the structures of $\sigma_{n-1}$ inspired by recent works of Guan, Ren and Wang \cite{GRW} and Ren and Wang \cite{RW}. It enables us to take care of the third order terms and therefore obtain the desired estimates. 

\medskip

The organization of the paper is as follows. In Section 2, we collect some basic properties for the asymptotic Plateau problem and an important inequality by Ren and Wang \cite{RW}. In Section 3, we prove Theorem \ref{theo-estimate} and Theorem \ref{theo-exist}.

\section{Preliminaries}

In this section, we first collect some basic properties for hypersurfaces in $\mathbb{H}^{n+1}$. We will use upper half-space model of $\mathbb{H}^{n+1}$.

Let $\Sigma$ be a connected, orientable, complete hypersurface in $\mathbb{H}^{n+1}$ with compact asymptotic boundary at infinity. Let $\mathbf{n}$ be the unit normal vector of $\Sigma$ pointing to the unbounded region in $\mathbb{R}^{n+1}_+\setminus \Sigma$.

Let $X$ and $\nu$ be the position vector and Euclidean unit normal vector of $\Sigma$ in $\mathbb{R}^{n+1}$, define
\begin{align*}
u=X\cdot \mathbf{e},\quad \nu^{n+1}=\mathbf{e}\cdot\nu,
\end{align*}
where $\mathbf{e}$ is the unit vector field in the positive $x_{n+1}$ direction in $\mathbb{R}^{n+1}$ and $\cdot$ denotes the Euclidean inner product in $\mathbb{R}^{n+1}$.

Let $e_1,\cdots,e_n$ be a local frame, then the metric and second fundamental form of $\Sigma$ are given by
\begin{align*}
g_{ij}=\left\langle e_i,e_j\right\rangle,\quad h_{ij}=\left\langle \overline{D} _{e_i}e_j,\mathbf{n}\right\rangle,
\end{align*}
where $\overline{D}$ denotes the Levi-Civita connection of $\mathbb{H}^{n+1}$.

In the following, we will assume $e_1,\cdots,e_n$ are orthonormal. Consequently, $g_{ij}=\delta_{ij}$ and $h_{ij}=\kappa_i\delta_{ij}$, where $\kappa_1,\cdots,\kappa_n$ are the principal curvatures of $\Sigma$.

The Gauss and Codazzi equations are given by
\begin{align*}
R_{ijkl}=-(\delta_{ik}\delta_{jl}-\delta_{il}\delta_{jk})+h_{ik}h_{jl}-h_{il}h_{jk}, 
\end{align*}
\begin{align*}
h_{ijk}=h_{ikj}.
\end{align*}
The convention that $R_{ijij}$ denotes the sectional curvature is used here.

The commutator formulas are given by
\begin{align}\label{comm}
h_{klij}=&\ h_{ijkl}-h_{ml}(h_{im}h_{kj}-h_{ij}h_{mk})-h_{mj}(h_{mi}h_{kl}-h_{il}h_{mk})\\ \nonumber
&-h_{ml}(\delta_{ij}\delta_{km}-\delta_{ik}\delta_{jm}) -h_{mj}(\delta_{il}\delta_{km}-\delta_{ik}\delta_{lm}).
\end{align}

\medskip

Denote
\begin{align*}
\sigma_{k}^{ii}=\frac{\partial \sigma_k}{\partial \kappa_i},\quad \sigma_{k}^{ii,jj}=\frac{\partial^2 \sigma_{k}}{\partial \kappa_i\partial\kappa_j}.
\end{align*}

We now collect two lemmas relating to $\nu^{n+1}$.
\begin{lemm}\label{nu}
\begin{align*}
\sum_i \frac{u_i^2}{u^2}=1-(\nu^{n+1})^2\leq 1,\quad  (\nu^{n+1})_i=\frac{u_i}{u}(\nu^{n+1}-\kappa_i),
\end{align*}
\begin{align*}
\sum_i \sigma_{n-1}^{ii}(\nu^{n+1})_{ii} =&\ 2\sum_i \sigma_{n-1}^{ii}\frac{u_i}{u} (\nu^{n+1})_i+(n-1)\sigma (1+(\nu^{n+1})^2)\\
&-\nu^{n+1}\left(\sum_i \sigma_{n-1}^{ii}+\sum_i \sigma_{n-1}^{ii}\kappa_i^2\right).
\end{align*}
\end{lemm}

\begin{proof}
The first one follows from the definition of $u$ and $\nu_{n+1}$. The second one can be found in the proof of Theorem 3.1 in \cite{GSX}. The third one can be found in (3.9) in \cite{GSX}. Note that our equation is not normalized as in \cite{GSX}, thus the format is slightly different from (3.9) in \cite{GSX}.
\end{proof}

\begin{lemm}\label{nua}
For $n\geq 2$, let $\Gamma=\partial\Omega$, where $\Omega$ is a bounded smooth domain in $\mathbb{R}^n$ with nonnegative mean curvature and $\sigma\in (0,n)$. Suppose $\Sigma$ is a $C^4$ vertical graph over $\Omega$ in $\mathbb{H}^{n+1}$ satisfying (\ref{sigma_n-1}). Then we have
\begin{align*}
\nu^{n+1}> a>0,
\end{align*}
where $a$ is a constant depending only on $n,\Omega$ and $\sigma$.
\end{lemm}

\begin{proof}
By Proposition 4.1 in \cite{GS10}, we have $w=\sqrt{1+|Du|^2}\leq C$. Since $\nu^{n+1}=\frac{1}{w}$, it follows that $\nu^{n+1}\geq \frac{1}{C}$.
\end{proof}

We now state some properties relating to $\sigma_k$ function.
\begin{lemm}\label{kn}
Assume $\kappa\in \Gamma_k$ and $\kappa_1\geq \cdots\geq \kappa_n$. Then
\begin{enumerate}
\item \begin{align*}
\sum_i \sigma_{k}^{ii}\kappa_i^2\geq \frac{k}{n}\sigma_1\sigma_k.
\end{align*}
\item If $\kappa_i\leq 0$, then 
\begin{align*}
-\kappa_i\leq \frac{n-k}{k}\kappa_1.
\end{align*}
\end{enumerate}

\end{lemm}

\begin{proof}
For (1), by Proposition 2.2 in \cite{HS}, we have
\begin{align*}
\sum_i \sigma_{k}^{ii}\kappa_i^2=\sigma_1\sigma_k-(k+1)\sigma_{k+1}.
\end{align*}
The inequality now follows by a Newton-Maclaurin inequality.

For (2), see Lemma 10 in \cite{RW2}.
\end{proof}

We now state an important inequality by Ren and Wang (Theorem 11 in \cite{RW}), which is crucial to our estimate.
\begin{lemm}\label{renwang}
Assume $\kappa\in \Gamma_{n-1}$ and $\kappa_1\geq \cdots\geq \kappa_n$. For any constant $\epsilon>0$, there exists a large constant $K$ depending only on $\epsilon$ such that
\begin{align*}
\kappa_1\left(K(\sigma_{n-1})_1^2-\sum_{p\neq q}\sigma_{n-1}^{pp,qq}h_{pp1}h_{qq1}\right)+(1+\epsilon)\sum_{i\neq 1} \sigma_{n-1}^{ii}h_{ii1}^2-\sigma_{n-1}^{11}h_{111}^2\geq 0.
\end{align*}

\end{lemm}

\section{Proof of main theorems}

\textit{Proof of Theorem \ref{theo-estimate}.}

\begin{proof}

Consider the quantity 
\begin{align*}
Q=\ln \kappa_1-N\ln \nu^{n+1},
\end{align*}
where $\kappa_1$ is the largest principle curvature and $N$ is a large constant to be determined later.

Suppose $Q$ attains maximum at an interior point $X_0$. If $\kappa_1$ has multiplicity more than $1$, then $Q$ is not smooth at $X_0$. To overcome this difficulty, we apply a standard perturbation argument, see for instance \cite{Chu}. Let $g$ be the first fundamental form of $\Sigma$ and $D$ be the corresponding Levi-Civita connection. Choose an orthonormal frame $e_1,\cdots,e_n$ near $X_0$ such that at $X_0$, we have
\begin{align*}
D_{e_i}e_j=0,\quad h_{ij}=\delta_{ij}\kappa_i,\quad \kappa_1\geq \cdots\geq \kappa_n.
\end{align*}

Near $X_0$, define a new tensor $B$ by
\begin{align*}
B(V_1,V_2)=g(V_1,V_2)-g(V_1,e_1)g(V_2,e_1),
\end{align*}
for tangent vectors $V_1$ and $V_2$. 

Denote $B_{ij}=B(e_i,e_j)$ and define
\begin{align*}
\tilde{h}_{ij}=h_{ij}-B_{ij}.
\end{align*}

Let $\tilde{\kappa}_1\geq \cdots\geq \tilde{\kappa}_n$ be the corresponding eigenvalues of $\tilde{h}_{ij}$.

It follows that $\kappa_1\geq \tilde{\kappa}_1$ near $X_0$ and at $X_0$, we have
\begin{align*}
\tilde{\kappa}_i=\begin{cases}
\kappa_1,\quad &i=1,\\
\kappa_i-1, & i>1.
\end{cases}
\end{align*}

Now consider the new test function 
\begin{align*}
\tilde{Q}=\ln\tilde{\kappa}_1-N\ln \nu^{n+1}.
\end{align*}

It also attains maximum at $X_0$. Moreover, $\tilde{\kappa}_1$ has multiplicity $1$, thus $\tilde{Q}$ is smooth at $X_0$.

At $X_0$, we have
\begin{align}\label{critical}
0=\frac{(\tilde{\kappa}_1)_i}{\tilde{\kappa}_1}-N\frac{(\nu^{n+1})_i}{\nu^{n+1}},
\end{align}
\begin{align}\label{first ineq}
0\geq \frac{(\tilde{\kappa}_1)_{ii}}{\tilde{\kappa}_1}-\frac{(\tilde{\kappa}_1)_i^2}{\tilde{\kappa}_1^2}-N\frac{(\nu^{n+1})_{ii}}{\nu^{n+1}}+N\frac{(\nu^{n+1})_i^2}{(\nu^{n+1})^2}.
\end{align}

By the definition of $B$ and the fact that $D_{e_i}e_j=0$ at $X_0$, we have
\begin{align*}
(\tilde{\kappa}_1)_i=\tilde{h}_{11i}=h_{11i},
\end{align*}
\begin{align*}
(\tilde{\kappa}_1)_{ii}=\tilde{h}_{11ii}+2\sum_{p\neq 1}\frac{\tilde{h}_{1pi}^2}{\tilde{\kappa}_1-\tilde{\kappa}_p}=h_{11ii}+2\sum_{p\neq 1}\frac{h_{1pi}^2}{\kappa_1-\tilde{\kappa}_p}.
\end{align*}

Plug into (\ref{first ineq}), we have
\begin{align*}
0\geq  \frac{h_{11ii}}{\kappa_1}+2\sum_{p\neq 1}\frac{h_{1pi}^2}{\kappa_1 (\kappa_1-\tilde{\kappa}_p)}-\frac{ h_{11i}^2}{\kappa_1^2}-N\frac{(\nu^{n+1})_{ii}}{\nu^{n+1}}.
\end{align*}

By commutator formula (\ref{comm}), we have
\begin{align*}
h_{11ii}=h_{ii11}+\kappa_1 ^2\kappa_i -\kappa_1\kappa_i ^2-\kappa_1+\kappa_i.
\end{align*}

It follows that
\begin{align*}
0\geq &\ \frac{h_{ii11}}{\kappa_1}+2\sum_{p\neq 1}\frac{h_{1pi}^2}{\kappa_1 (\kappa_1-\tilde{\kappa}_p)}-\frac{ h_{11i}^2}{\kappa_1^2}-N\frac{(\nu^{n+1})_{ii}}{\nu^{n+1}}\\
&+\kappa_1 \kappa_i -\kappa_i ^2-1+\frac{\kappa_i}{\kappa_1}.
\end{align*}

Contract with $F^{ii}=\sigma_{n-1}^{ii}$, we have
\begin{align}\label{second ineq}
0\geq &\sum_i \frac{F^{ii}h_{ii11}}{\kappa_1}+2\sum_i\sum_{p\neq 1}\frac{F^{ii}h_{1pi}^2}{\kappa_1 (\kappa_1-\tilde{\kappa}_p)}-\sum_i\frac{ F^{ii}h_{11i}^2}{\kappa_1^2}\\\nonumber
&-N\sum_i\frac{F^{ii}(\nu^{n+1})_{ii}}{\nu^{n+1}}-\sum_i F^{ii}\kappa_i ^2-\sum_i F^{ii},
\end{align}
we have used the fact $\sum_i F^{ii}\kappa_i=(n-1)\sigma>0$ in the above inequality.

By Lemma \ref{nu}, we have
\begin{align*}
\sum_i F^{ii}(\nu^{n+1})_{ii} =&\ 2\sum_i F^{ii}\frac{u_i}{u} (\nu^{n+1})_i+(n-1)\sigma (1+(\nu^{n+1})^2)\\
&-\nu^{n+1}\left(\sum_i F^{ii}+\sum_i F^{ii}\kappa_i^2\right).
\end{align*}

Plug into (\ref{second ineq}), we have
\begin{align*}
0\geq &\sum_i \frac{F^{ii}h_{ii11}}{\kappa_1}+2\sum_i\sum_{p\neq 1}\frac{F^{ii}h_{1pi}^2}{\kappa_1 (\kappa_1-\tilde{\kappa}_p)}-\sum_i\frac{ F^{ii}h_{11i}^2}{\kappa_1^2}-2N \sum_i  F^{ii}\frac{u_i}{u} \frac{(\nu^{n+1})_i}{\nu^{n+1}} \\
&+ (N-1)\left(\sum_i F^{ii}+\sum_i F^{ii}\kappa_i^2\right)-N(n-1)\sigma \frac{1+(\nu^{n+1})^2}{\nu^{n+1}}.
\end{align*}

By Lemma \ref{nua}, we have
\begin{align}\label{third ineq}
0\geq &\sum_i \frac{F^{ii}h_{ii11}}{\kappa_1}+2\sum_i\sum_{p\neq 1}\frac{F^{ii}h_{1pi}^2}{\kappa_1 (\kappa_1-\tilde{\kappa}_p)}-\sum_i\frac{ F^{ii}h_{11i}^2}{\kappa_1^2}\\\nonumber
&-2N \sum_i  F^{ii}\frac{u_i}{u} \frac{(\nu^{n+1})_i}{\nu^{n+1}}+ (N-1)\left(\sum_i F^{ii}+\sum_i F^{ii}\kappa_i^2\right)-CN,
\end{align}
where $C$ is a universal constant depending only on $n, \Omega$ and $\sigma$. From now on, we will use $C$ to denote a universal constant depending only on $n, \Omega$ and $\sigma$, it may change from line to line.

Differentiate (\ref{sigma_n-1}), we have
\begin{align*}
\sum_i F^{ii}h_{ii11}=-\sum_{p,q,r,s} F^{pq,rs}h_{pq1}h_{rs1}=-\sum_{p\neq q}F^{pp,qq}h_{pp1}h_{qq1}+\sum_{p\neq q}F^{pp,qq}h_{pq1}^2.
\end{align*}

Plug into (\ref{third ineq}), we have
\begin{align}\label{fourth ineq}
0\geq& -\sum_{p\neq q}\frac{F^{pp,qq}h_{pp1}h_{qq1}}{\kappa_1}+\sum_{p\neq q} \frac{F^{pp,qq}h_{pq1}^2}{\kappa_1 }+2\sum_i\sum_{p\neq 1}\frac{F^{ii}h_{1pi}^2}{\kappa_1 (\kappa_1-\tilde{\kappa}_p)}\\\nonumber
&-\sum_i\frac{ F^{ii}h_{11i}^2}{\kappa_1^2}-2N \sum_i  F^{ii}\frac{u_i}{u} \frac{(\nu^{n+1})_i}{\nu^{n+1}}\\\nonumber
&+(N-1)\left(\sum_i F^{ii}+\sum_i F^{ii}\kappa_i^2\right)-CN.
\end{align}

Now
\begin{align*}
2\sum_i\sum_{p\neq 1}\frac{F^{ii}h_{1pi}^2}{\kappa_1 (\kappa_1-\tilde{\kappa}_p)}&\geq 2\sum_{p\neq 1}\frac{F^{pp}h_{1pp}^2}{\kappa_1 (\kappa_1-\tilde{\kappa}_p)}+2\sum_{p\neq 1}\frac{F^{11}h_{1p1}^2}{\kappa_1 (\kappa_1-\tilde{\kappa}_p)}\\
&=2\sum_{i\neq 1}\frac{F^{ii}h_{ii1}^2}{\kappa_1 (\kappa_1-\tilde{\kappa}_i)}+2\sum_{i\neq 1}\frac{F^{11}h_{11i}^2}{\kappa_1 (\kappa_1-\tilde{\kappa}_i)}.
\end{align*}

Without loss of generality, assume $\kappa_1$ has multiplicity $m$, then we have
\begin{align*}
\sum_{p\neq q} \frac{F^{pp,qq}h_{pq1}^2}{\kappa_1 }&\geq 2\sum_{i>m}\frac{F^{11,ii}h_{11i}^2}{\kappa_1 }=2\sum_{i>m}\frac{(F^{ii}-F^{11})h_{11i}^2}{\kappa_1 (\kappa_1 -\kappa_i )}\geq 2\sum_{i>m}\frac{(F^{ii}-F^{11})h_{11i}^2}{\kappa_1 (\kappa_1 -\tilde{\kappa}_i )},
\end{align*}
where we have used the fact $\tilde{\kappa}_i=\kappa_i-1$ in the last inequality.

Plug the above two inequalities into (\ref{fourth ineq}), we have
\begin{align}\label{fifth ineq}
0\geq& -\sum_{p\neq q}\frac{F^{pp,qq}h_{pp1}h_{qq1}}{\kappa_1}+2\sum_{i\neq 1}\frac{F^{ii}h_{ii1}^2}{\kappa_1 (\kappa_1-\tilde{\kappa}_i)}+2\sum_{1<i\leq  m}\frac{F^{11}h_{11i}^2}{\kappa_1 (\kappa_1-\tilde{\kappa}_i)}\\\nonumber
&+2\sum_{i>m}\frac{F^{ii}h_{11i}^2}{\kappa_1 (\kappa_1 -\tilde{\kappa}_i )}-\sum_i\frac{ F^{ii}h_{11i}^2}{\kappa_1^2}-2N \sum_i  F^{ii}\frac{u_i}{u} \frac{(\nu^{n+1})_i}{\nu^{n+1}}\\\nonumber
&+ (N-1)\left(\sum_i F^{ii}+\sum_i F^{ii}\kappa_i^2\right)-CN.
\end{align}

In the following, we will always assume $\kappa_1$ is sufficiently large, for otherwise we have obtained the estimate.

For $1<i\leq m$, we have
\begin{align*}
2\frac{F^{11}h_{11i}^2}{\kappa_1 (\kappa_1-\tilde{\kappa}_i)}-\frac{ F^{ii}h_{11i}^2}{\kappa_1^2}=(2\kappa_1-1) \frac{ F^{ii}h_{11i}^2}{\kappa_1^2}\geq \frac{ F^{ii}h_{11i}^2}{\kappa_1^2}.
\end{align*}

For $i>m$, by Lemma \ref{kn}, we have
\begin{align*}
2\frac{F^{ii}h_{11i}^2}{\kappa_1 (\kappa_1 -\tilde{\kappa}_i )}-\frac{ F^{ii}h_{11i}^2}{\kappa_1^2}=\frac{F^{ii} (\kappa_1 +\tilde{\kappa}_i )}{\kappa_1^2 (\kappa_1 -\tilde{\kappa}_i )}h_{11i}^2\geq c(n)\frac{ F^{ii}h_{11i}^2}{\kappa_1^2},
\end{align*}
where $c(n)$ is a constant depending only on $n$.

Plug into (\ref{fifth ineq}), we have
\begin{align*}
0\geq  &-\sum_{p\neq q}\frac{F^{pp,qq}h_{pp1}h_{qq1}}{\kappa_1}+2\sum_{i\neq 1}\frac{F^{ii}h_{ii1}^2}{\kappa_1 (\kappa_1-\tilde{\kappa}_i)}+ c(n)\sum_{i\neq 1}\frac{F^{ii}h_{11i}^2}{\kappa_1^2}-\frac{ F^{11}h_{111}^2}{\kappa_1^2}\\
&-2N \sum_i  F^{ii}\frac{u_i}{u} \frac{(\nu^{n+1})_i}{\nu^{n+1}}+ (N-1)\left(\sum_i F^{ii}+\sum_i F^{ii}\kappa_i^2\right)-CN.
\end{align*}

By Lemma \ref{kn}, the fact that $n\geq 3$ and that $\kappa_1$ is sufficiently large, we have
\begin{align*}
\frac{2}{\kappa_1(\kappa_1-\tilde{\kappa}_i)}\geq \frac{2}{\kappa_1\left(\frac{n}{n-1}\kappa_1+1\right) }\geq \frac{7}{6}\frac{1}{\kappa_1^2}.
\end{align*}

Apply Lemma \ref{renwang}, we have
\begin{align*}
0\geq  c(n)\sum_{i\neq 1}\frac{F^{ii}h_{11i}^2}{\kappa_1^2}-2N \sum_i  F^{ii}\frac{u_i}{u} \frac{(\nu^{n+1})_i}{\nu^{n+1}}+(N-1)\left(\sum_i F^{ii}+\sum_i F^{ii}\kappa_i^2\right)-CN.
\end{align*}

By Lemma \ref{nu}, we have
\begin{align*}
-2N F^{ii}\frac{u_i}{u} \frac{(\nu^{n+1})_i}{\nu^{n+1}}=-2N F^{ii}\frac{u_i^2}{u^2} \frac{\nu^{n+1}-\kappa_i}{\nu^{n+1}}.
\end{align*}

It is negative only if $\kappa_i<\nu^{n+1}$. Since $\kappa_1$ is sufficiently large, in particular $\kappa_1>\nu^{n+1}$, it follows that
\begin{align*}
0\geq  c(n)\sum_{i\neq 1}\frac{F^{ii}h_{11i}^2}{\kappa_1^2}-2N \sum_{i\neq 1} F^{ii}\frac{u_i}{u} \frac{(\nu^{n+1})_i}{\nu^{n+1}}+(N-1)\left(\sum_i F^{ii}+\sum_i F^{ii}\kappa_i^2\right)-CN.
\end{align*}

By critical equation (\ref{critical}), we have
\begin{align*}
0\geq &\  c(n)N^2\sum_{i\neq 1} F^{ii}\frac{(\nu^{n+1})_i^2}{(\nu^{n+1})^2}-2N \sum_{i\neq 1}  F^{ii}\frac{u_i}{u} \frac{(\nu^{n+1})_i}{\nu^{n+1}}\\
&+ (N-1)\left(\sum_i F^{ii}+\sum_i F^{ii}\kappa_i^2\right)-CN\\
\geq & -C\sum_{i\neq 1} F^{ii}\frac{u_i^2}{u^2}+(N-1)\left(\sum_i F^{ii}+\sum_i F^{ii}\kappa_i^2\right)-CN.
\end{align*}

Choose $N$ sufficiently large, together with Lemma \ref{nu} and Lemma \ref{kn}, we have
\begin{align*}
0 \geq  (N-1)\sum_i F^{ii}\kappa_i^2-CN \geq (N-1) \frac{n-1}{n}\sigma\sigma_1-CN.
\end{align*}

It follows that $\sigma_1\leq C$. The theorem is now proved.

\end{proof}

\textit{Proof of Theorem \ref{theo-exist}.}

\begin{proof}
As pointed out in \cite{GS10}, the only missing piece is the interior curvature estimates. In view of Theorem \ref{theo-estimate}, Theorem \ref{theo-exist} is now proved, see details in \cite{GS10}.
\end{proof}

\end{document}